\documentclass[12pt]{article}

\newtheorem{theorem}{Theorem}

\newtheorem{corollary}[theorem]{Corollary}

\newtheorem{remark}[theorem]{Remark}

\newenvironment{proof}[1][Proof]{\textbf{#1.} }{\ \rule{0.5em}{0.5em}}
\newdimen\dummy
\dummy=\oddsidemargin
\begin{document}

\title{The Hawaiian earring group and metrizability}
\author{Paul Fabel \\
Department of Mathematics \& Statistics\\
Mississippi State University}
\date{}
\maketitle

\begin{abstract}
Endowed with quotient topology inherited from the space of based loops, the
fundamental group of the Hawaiian earring fails to be metrizable.  The
fundamental group of any space which retracts to the Hawaiian earring is
also nonmetrizable.
\end{abstract}

\section{Introduction}

The Hawaiian earring $HE$ is the union of a null sequence of simple closed
curves joined at a common point $p,$ and is among the simplest of spaces to
which the standard tools of algebraic topology fail to apply.

The fundamental group $\pi _{1}(HE)$ has received attention from a various
algebraic and combinatorial perspectives (\cite{can} \cite{desmit}\cite{eda3}%
). This paper settles a fundamental question regarding the \textbf{topology}
of $\pi _{1}(HE).$

There are at least 3 natural ways of imparting a topology on $\pi _{1}(HE).$

\begin{enumerate}
\item  It is a (nontrivial) fact that the natural homomorphism into the
inverse limit of (discrete) free groups $\phi :\pi _{1}(HE)\rightarrow
\lim_{\leftarrow }\pi _{1}(X_{n})$ \ is one to one \cite{mor}, and thus one
can pull back via the embedding $\phi $ to create the metric space $\pi
_{1}^{\lim }(HE,p)$

\item  One can consider two loop classes $[f]$ and $\ [g]$ to be close (\cite
{bog} \cite{FZ}) if there is a small path $\alpha $ such that $f\ast \alpha $
is path homotopic to $g$ . This creates a metric space $\pi _{1}^{met}(HE,p).
$

\item  One can endow the path components of the loop space over $HE$ with
the quotient topology. Denote this space $\pi _{1}^{top}(HE,p)$ as in \cite
{Biss}.
\end{enumerate}

It is not obvious that these determine distinct spaces. Indeed \cite{fab3}
is devoted to a proof that $\pi _{1}^{\lim }(HE,p)$ and $\pi
_{1}^{top}(HE,p) $ are not \textbf{canonically} homeomorphic.

It is shown in \cite{fab1} that $\pi _{1}^{top}(HE,p)$ is separable and
normal, but the question of metrizability was not resolved. The main result
of this paper is that $\pi _{1}^{top}(HE)$ is \textbf{not} metrizable.

Thus we can conclude that $\pi _{1}^{top}(HE,p)$ is topologically distinct
from the other two spaces. Moreover it follows from a ``no retraction''
theorem of the author \cite{fab2} that $HE$ cannot be embedded as a retract
of any space $Y$ such that $\pi _{1}^{top}(Y)$ is metrizable. For example if 
$Y$ is the Sierpinski Gasket or the Menger Sponge, the topological
fundamental group $\pi _{1}^{top}(Y)$ fails to be metrizable.

\section{Main Result}

Suppose $X$ is a metric space and $p\in X.$ Let $C_{p}(X)=\{f:[0,1]%
\rightarrow X$ such that $f$ is continuous and $f(0)=f(1)=p\}.$ Endow $%
C_{p}(X)$ with the topology of uniform convergence. Let $P$ denote the
constant path and let $[P]$ denote the path component of $P$ in $C_{p}(X).$

The \textbf{topological fundamental group} $\pi _{1}^{top}(X,p)$ is the set
of path components of $C_{p}(X)$ endowed with the quotient topology under
the canonical surjection $\Pi :C_{p}(X)\rightarrow \pi _{1}^{top}(X,p)$
satisfying $\Pi (f)=\Pi (g)$ if and only if $f$ and $g$ belong to the same
path component of $C_{p}(X).$

Thus a set $U\subset \pi _{1}^{top}(X,p)$ is open in $\pi _{1}^{top}(X,p)$
if and only if $q^{-1}(U)$ is open in $C_{p}(X).$

Let $X_{n}$ denote the circle of radius $\frac{1}{n}$ centered at $(\frac{1}{%
n},0)$ in the plane.

Let $HE=\cup _{n=1}^{\infty }X_{n},$ denote the familiar Hawaiian earring,
with $p=(0,0).$ Let $f\ast g$ denote the familiar operation of path
concatenation.

\begin{remark}
The topological Hawaiian earring group $\pi _{1}^{top}(HE,p)$ is separable
and normal. See \cite{fab1}.
\end{remark}

\begin{theorem}
\label{main}The topological Hawaiian earring group $\pi _{1}^{top}(HE,p)$ is
not metrizable.
\end{theorem}

\begin{proof}
The fact that $\pi _{1}^{top}(HE,p)$ is a $T_{1}$ space is equivalent (\cite
{fab1}) to the fact (\cite{mor}) that the natural homomorphism $\phi :\pi
_{1}^{top}(HE,p)\rightarrow \lim_{\leftarrow }\pi _{1}(X_{n},p)$ is one to
one. 

Suppose $U_{1},U_{2},,.$ is a sequence of open sets in $\pi _{1}^{top}(X,p)$
such that $[P]\in U_{n}$ for each $n$. Let $V_{n}=\Pi ^{-1}(U_{n}).$ Because 
$P\in V_{n},$ the path homotopic map $(\alpha _{1}\ast P\ast \alpha
_{1}^{-1}\ast P)^{n}\in V_{n}.$  Because $V_{n}$ is open in $C_{p}(HE)$
there exists $m_{n}>1$ such that $(\alpha _{1}\ast \alpha _{m_{n}}\ast
\alpha _{1}^{-1}\ast \alpha _{m_{n}}^{-1})^{n}\in V_{n}.$ (In this case $%
\alpha _{i}$ and $\alpha _{i}^{-1}$ denote loops that orbit once around $%
X_{i}$ in opposite directions.)

Let $f_{n}=(\alpha _{1}\ast \alpha _{m_{n}}\ast \alpha _{1}^{-1}\ast \alpha
_{m_{n}}^{-1})^{n}.$ If we can show that the sequence $[f_{n}]$ does not
converge to $[P]$ in $\pi _{1}^{top}(HE,p),$ it will follow immediately that 
$U_{n}$ cannot be the open metric ball of radius $\frac{1}{n}$ centered at $%
[P],$ and in particular $\pi _{1}^{top}(HE,p)$ cannot be metrizable.

The proof that $\{[f_{n}]\}$ does not converge in $\pi _{1}^{top}(HE,p)$ to $%
[P]$ is very similar to the proof \cite{fab3} that $\pi _{1}^{top}(HE,p)$
does not naturally embed topologically in $\lim_{\leftarrow }\pi
_{1}(X_{n},p)$. Choose $q\neq p$ such that $q\in X_{1}.$ Define an
oscillation function $o:C_{p}(HE)\rightarrow \{1,2,3,...\}$ such $o(g)=M$ if 
$M$ is the largest nonnegative integer such that there exists $%
t_{0}<t_{1}<...t_{M}$ with $g(t_{2i})=p$ and $g(t_{2i+1})=q.$

Notice if $g_{n}$ is path homotopic to $f_{n}$ in $HE,$ then $o(g_{n})\geq
o(f_{n})\geq n.$ (To see this, first replace $g_{n}$ by $r_{n}(g_{n})$ where 
$r_{n}:HE\rightarrow X_{n}$ is the natural retraction and note $%
o(r_{n}g_{n})\leq o(g_{n})$ and moreover $r_{n}g_{n}$ is also path homotopic
to $f_{n}=r_{n}(f_{n}).$ Next observe that at each step (via the van Kampen
Theorem \cite{Hatch}) of the reduction of $r_{n}(g_{n})$ to $f_{n},$ the
oscillation number does not increase).

Next observe that if $g_{n}\rightarrow g$ in $C_{p}(HE)$ that $o(g)\geq
o(g_{n})$ eventually.

Thus, if $g_{n}$ is path homotopic to $f_{n}$ for each $n$ then the sequence 
$g_{1},g_{2},...$ cannot have a subsequential limit $g\in C_{p}(HE)$ since
it would follow that $o(g)=\infty .$ Moreover, since $\pi _{1}^{top}(HE,p)$
is a $T_{1}$ space, for each $n$ the path component of $f_{n}$ is closed in $%
C_{p}(HE).$

Hence the preimage under $\Pi $ of the sequence $[f_{1}],[f_{2}],..$
determines a closed subset of $C_{p}(HE),$ and since $\Pi $ is a quotient
map the set $\{[f_{1}],[f_{2}],..\}$ is closed in $\pi _{1}(HE,p)$. Notice $%
[P]$ is not in this closed set, and consequently $[P]$ is not a limit point
of this set, and hence the sequence $\{[f_{n}]\}$ does not converge to $[P]$
in $\pi _{1}^{top}(HE,p).$
\end{proof}

\begin{remark}
The paper \cite{fab3} shows that $\pi _{1}^{top}(HE,p)$ does not embed
naturally in the (metrizable) inverse limit space $\lim_{\leftarrow }\pi
_{1}(X_{n},p).$ Theorem \ref{main} reveals there exists no topological
embedding whatsoever and in particular $\pi _{1}^{top}(HE,p)$ and $\pi
_{1}^{\lim }(HE,p)$ are not homeomorphic.
\end{remark}

\begin{corollary}
\label{main2}Suppose $Y$ is a path connected metric space such that $\pi
_{1}^{top}(Y)$ is metrizable . Then the Hawaiian earring cannot be embedded
as a retract of $Y.$
\end{corollary}

\begin{proof}
An embedding of $HE$ into $Y$ would induce a topological embedding of $\pi
_{1}^{top}(HE)$ into $\pi _{1}^{top}(Y)$ as shown in \cite{fab2},
contradicting the nonmetrizability of $\pi _{1}(HE).$
\end{proof}


\begin{thebibliography}{99}
\bibitem{Biss}  Biss, Daniel K. The topological fundamental group and
generalized covering spaces. Topology Appl. 124 (2002), no. 3, 355--371.

\bibitem{bog}  Bogley, W.A., Sieradski, A.J. Universal path spaces. Preprint.

\bibitem{can}  Cannon, J. W.; Conner, G. R. The combinatorial structure of
the Hawaiian earring group. Topology Appl. 106 (2000), no. 3, 225--271.

\bibitem{desmit}  de Smit, Bart. The fundamental group of the Hawaiian
earring is not free. Internat. J. Algebra Comput. 2 (1992), no. 1, 33--37.

\bibitem{eda3}  Eda, Katsuya. Free subgroups of the fundamental group of the
Hawaiian earring. J. Algebra 219 (1999), no. 2, 598--605.

\bibitem{fab1}  Fabel, Paul A monomorphism theorem for the inverse limit of
nested retracts. Preprint. http://front.math.ucdavis.edu/math.AT/0502275

\bibitem{fab2}  Fabel, Paul A retraction theorem for topological fundamental
groups with application to the Hawaiian earring. Preprint
http://front.math.ucdavis.edu/math.AT/0502218

\bibitem{fab3}  Fabel, Paul. The topological Hawaiian earring group does not
embed in the inverse limit of free groups. Algebraic and Geometric Topology
5 (2005). http://www.maths.warwick.ac.uk/agt/AGTVol5/agt-5-64.abs.html

\bibitem{FZ}  Fischer, Hanspeter; Zastrow, Andreas. Generalized universal
coverings and the shape group. Preprint

\bibitem{Hatch}  Hatcher, Allen. Algebraic topology. Cambridge University
Press, Cambridge, 2002.

\bibitem{mor}  Morgan, John W.; Morrison, Ian A van Kampen theorem for weak
joins. Proc. London Math. Soc. (3) 53 (1986), no. 3, 562--576.
\end{thebibliography}
\end{document}